\documentclass[12pt,thmsa]{article}
\usepackage{amsfonts}
\usepackage{amsmath}
\usepackage{amssymb}
\usepackage{sw20lart}

\input{tcilatex}
\begin{document}

\author{Krzysztof Ma\'{s}lanka\thanks{%
e-mail:kmaslank@uj.edu.pl} \\
Polish Academy of Sciences\\
Institute of the History of Science\\
ul. \'{s}w. Jana 22, 31-018 Krak\'{o}w, Poland}
\title{B\'{a}ez-Duarte's Criterion for the Riemann Hypothesis and Rice's
Integrals}
\maketitle

\begin{abstract}
Criterion for the Riemann hypothesis found by B\'{a}ez-Duarte involves
certain real coefficients $c_{k\text{ }}$defined as alternating binomial
sums. These coefficients can be effectively investigated using N\"{o}%
rlund-Rice's integrals. Their behavior exhibits characteristic trend, due to
trivial zeros of zeta, and fading oscillations, due to complex zeros. This
method enables to calculate numerical values of $c_{k\text{ }}$for large
values of $k$, at least to $k=4\cdot 10^{8}$.

We give explicit expressions both for the trend and for the oscillations.
The first tends to zero and is therefore, in view of the criterion,
irrelevant for the Riemann hypothesis. The oscillations can be further
decomposed into a series of harmonics with amplitudes diminishing quickly.
Possible violation of the Riemann hypothesis would indicate that the
amplitude of some high harmonic increases.\newpage
\end{abstract}

\section{Introduction}

Several years ago a new expansion for the Riemann zeta function has been
found by the author \cite{Maslanka}, \cite{Duarte0}:%
\begin{equation}
\zeta (s)=\frac{1}{s-1}\dsum\limits_{k=0}^{\infty }\frac{\Gamma \left( k+1-%
\frac{s}{2}\right) }{\Gamma \left( 1-\frac{s}{2}\right) }\frac{A_{k}}{k!}%
\equiv \frac{1}{s-1}\dsum\limits_{k=0}^{\infty }\left( 1-\frac{s}{2}\right)
_{k}\frac{A_{k}}{k!}  \label{zeta}
\end{equation}%
where%
\begin{equation}
A_{k}:=\dsum\limits_{j=0}^{k}(-1)^{j}\binom{k}{j}(2j-1)\zeta (2j+2)\equiv
\dsum\limits_{j=0}^{k}\binom{k}{j}\frac{\pi ^{2j+2}}{\left( 2\right)
_{j}\left( \frac{1}{2}\right) _{j}}B_{2j+2}  \label{ak}
\end{equation}%
and%
\begin{equation*}
\left( x\right) _{k}\equiv \frac{\Gamma \left( k+x\right) }{\Gamma \left(
x\right) }
\end{equation*}%
is the Pochhammer symbol (having such an unfortunate denotation and
sometimes called the rising factorial). Pochhammer symbols are in fact
polynomials in variable $x$ with integer coefficients equal to the Stirling
numbers of the first kind. Since $A_{k}$ tend to zero fast enough as $%
k\rightarrow \infty $ expansion (\ref{zeta}) converges uniformly on the
whole complex plane.

In fact, there exists a whole class of expansions similar to (\ref{zeta}).
The crucial thing is to remove the single pole of $\zeta (s)$ in $s=1$ which
may be achieved either by multiplication by $s-1$ (or by any other function
having single simple zero at unity) or by subtraction of $\frac{1}{s-1}$ (or
by any other function having single simple pole at unity). In the second
case we get a series converging even faster:%
\begin{equation}
\zeta (s)=\frac{1}{s-1}+\dsum\limits_{k=0}^{\infty }\frac{\Gamma \left( k+1-%
\frac{s}{2}\right) }{\Gamma \left( 1-\frac{s}{2}\right) }\frac{A_{k}}{k!}
\label{zeta1}
\end{equation}%
where%
\begin{equation}
A_{k}:=\dsum\limits_{j=0}^{k}(-1)^{j}\binom{k}{j}\left( \zeta (2j+2)-\frac{1%
}{2j+1}\right)  \label{ak1}
\end{equation}

One has also the freedom of choosing the "node" points at which expressions (%
\ref{zeta}) or (\ref{zeta1}) give exact values, however the choice of even
positive integers seems most natural, since $A_{k}$ may be expressed by
Bernoulli numbers and appropriate powers of $\pi $ avoiding values such as $%
\zeta (3)$.

Finally, on can go all the way writing 
\begin{equation}
\zeta (s)=\frac{1}{2(s-1)}\left[ 1+s\left( \log (2\pi
)-1+s\dsum\limits_{k=0}^{\infty }\frac{\Gamma \left( k+1-\frac{s}{2}\right) 
}{\Gamma \left( 1-\frac{s}{2}\right) }\frac{A_{k}}{k!}\right) \right]
\label{zeta2}
\end{equation}%
with 
\begin{equation*}
A_{k}:=\dsum\limits_{j=0}^{k}(-1)^{j}\binom{k}{j}f(2j+2)
\end{equation*}%
and%
\begin{equation*}
f(s):=\frac{\frac{2(s-1)\zeta (s)-1}{s}+1-\log (2\pi )}{s}
\end{equation*}%
where $f(s)$ is also regular on the whole complex plane. One can check that
expansion (\ref{zeta2}) gives exact values $\zeta (0)=-\frac{1}{2}$ as well
as $\zeta ^{\prime }(0)=-\frac{1}{2}\log (2\pi )$ irrespective of how many
terms in the series is taken into account.

In 2003 Luis B\'{a}ez-Duarte, investigating expansion (\ref{zeta}), found an
interesting criterion for the Riemann hypothesis (RH) \cite{Duarte1}. The
crucial thing is to estimate the asymptotic behavior of certain real numbers
defined as%
\begin{equation}
c_{k}:=\sum_{j=0}^{k}(-1)^{j}\binom{k}{j}\frac{1}{\zeta (2j+2)}.  \label{ck}
\end{equation}%
More precisely, RH is equivalent to the statement that%
\begin{equation}
c_{k}\ll k^{-3/4+\varepsilon },\forall \varepsilon >0.  \label{criterion}
\end{equation}%
If 
\begin{equation*}
c_{k}\ll k^{-3/4}
\end{equation*}%
then all zeros of $\zeta $ are simple whereas regardless of the validity of
the RH (unconditionally)%
\begin{equation*}
c_{k}\ll k^{-1/2}.
\end{equation*}%
As B\'{a}ez-Duarte pointed out, criteria of this type, i.e. relating RH to
values of zeta at integer points, were known for a long time, however this
one is definitely simpler and well-fitted for numerical investigations. Such
investigations were performed numerically by the author in 2003 leading to
the observation (which was surprising for us at that time) that the global
behavior of coefficients $c_{k}$ may be split into a trend (which dominates
for low $k$ not exceeding roughly $10^{4}$) and subtle oscillations
superimposed on this trend. (Similar splitting has been found by the author
in the case of Li criterion for the RH, however in this case the
oscillations are much more chaotic \cite{Maslanka2}.)

\section{Binomial transforms and asymptotics}

Of course, such behavior may be well understood in the theory of binomial
transforms using Rice's integrals (already known by N\"{o}rlund). Sums of
this type can be interpreted as high order differences of appropriate
numerical sequences $\left\{ \varphi _{k}\right\} $. They were carefully
investigated \cite{Flajolet}, see also \cite{Preu}. The main result is (cf.
Theorem 2 in \cite{Flajolet}):%
\begin{equation}
\sum_{k=n_{0}}^{n}(-1)^{k}\binom{n}{k}\varphi (k)=-(-1)^{n}\sum_{s}\func{Res}%
\left[ \varphi (s)\frac{n!}{s(s-1)...(s-n)}\right]  \label{theorem}
\end{equation}

The proof of this fundamental theorem consists in applying classic Cauchy
residue formula.

In the case of (\ref{ck})\ we have $n_{0}=0$ and $\varphi (s)=1/\zeta (2s+2)$%
, which is analytic on $[0,\infty \lbrack $ and meromorphic on the whole
complex plane $%
\mathbb{C}
$, therefore it fulfils the assumptions of the above theorem.

\section{Results}

It is clear that the problem is to find all poles of suitable function, in
our case $1/\zeta (s)$. These poles are due to simple zeros of zeta in $%
s=-2n $ as well as due to complex zeros which we write as $\rho =\frac{1}{2}%
\pm i\gamma $. (Obviously if RH were true then all $\gamma $s would be
real.) It is a matter of elementary exercise to find that for $n=1,2,...$
the residues in real poles are:%
\begin{equation}
\func{Res}\left( \frac{1}{\zeta (s)};s=-2n\right) =\frac{1}{\zeta ^{\prime
}(-2n)}=2\frac{(-1)^{n}(2\pi )^{2n}}{(2n)!\zeta (2n+1)}  \label{trivialres}
\end{equation}%
where the last equality is a consequence of the functional equation for the
zeta function.

Using main theorem (\ref{theorem}) and (\ref{trivialres}) we find that the
asymptotic form of the trend which stems from trivial zeros is%
\begin{eqnarray}
\overset{\_}{c}_{k} &=&-\frac{k!}{\pi ^{3/2}}\sum_{m=2}^{\infty }\frac{%
(-1)^{m}\pi ^{2m}}{\Gamma (k+m+1)\Gamma (m-\frac{1}{2})}\frac{1}{\zeta (2m-1)%
}  \label{ctrend} \\
&=&-\frac{1}{4\pi ^{2}}\sum_{m=2}^{\infty }\frac{B(k+1,m)}{\Gamma (2m-1)}%
\frac{(-1)^{m}(2\pi )^{2m}}{\zeta (2m-1)}  \notag
\end{eqnarray}%
(the latter being numerically still more effective) whereas the oscillating
part due to complex zeros in the critical strip is:%
\begin{eqnarray}
\overset{\sim }{c}_{k} &=&k!\sum_{\rho }\frac{\Gamma \left( \frac{1+\rho }{2}%
\right) }{\Gamma \left( k+1+\frac{1+\rho }{2}\right) }\func{Res}\left( \frac{%
1}{\zeta (2s+2)};s=-\frac{1+\rho }{2}\right)  \label{osc} \\
&=&\sum_{\rho }B\left( k+1,\frac{1+\rho }{2}\right) \func{Res}\left( \frac{1%
}{\zeta (2s+2)};s=-\frac{1+\rho }{2}\right)  \notag
\end{eqnarray}%
where $B(x,y)$ is the Euler beta function and $c_{k}=\overset{\_}{c}_{k}+%
\overset{\sim }{c}$. Both (\ref{ctrend}) and (\ref{osc}) converge quickly,
hence they are suitable for numerical estimations contrary to the direct
approach using the main definition (\ref{ck}) which becomes very
time-consuming as $k$ grows. What's more, values of zeta function at
positive even integers should be calculated with many significant digits.
Since $\zeta (2n)$ tends quickly to unity as $n$ grows it is advisable, when
using (\ref{ck}), to tabulate appropriate number of high-precision values of 
$\zeta (2n)-1$ to preserve sufficient amount of significant digits as well
as to avoid repeated unnecessary calculations. Nevertheless, calculating of
the single coefficient $c_{300000}$ took about 2 weeks on a fast cluster of
four computers \cite{Wolf}. On the other hand, formulas (\ref{ctrend}) and (%
\ref{osc}) are much more effective: obtaining $c_{1000000}$ is a matter of
few tens of seconds on a modest machine. In numerical calculations using (%
\ref{osc}) the function \texttt{NResidue} implemented in \textit{Mathematica}
proved especially useful. Figures 1--5 and the table below present the
results.\bigskip

\begin{center}
\begin{tabular}{|c|c|}
\hline
$k$ & $c_{k}$ \\ \hline
$10^{5}$ & $+1.60976\cdot 10^{-9}$ \\ \hline
$2\cdot 10^{5}$ & $-7.89739\cdot 10^{-9}$ \\ \hline
$3\cdot 10^{5}$ & $+5.82876\cdot 10^{-9}$ \\ \hline
$4\cdot 10^{5}$ & $-2.89364\cdot 10^{-9}$ \\ \hline
$5\cdot 10^{5}$ & $-3.45567\cdot 10^{-9}$ \\ \hline
$6\cdot 10^{5}$ & $+1.13652\cdot 10^{-9}$ \\ \hline
$7\cdot 10^{5}$ & $+3.14429\cdot 10^{-9}$ \\ \hline
$8\cdot 10^{5}$ & $+2.00526\cdot 10^{-9}$ \\ \hline
$9\cdot 10^{5}$ & $-1.70316\cdot 10^{-10}$ \\ \hline
$10^{6}$ & $-1.77502\cdot 10^{-9}$ \\ \hline
$2\cdot 10^{6}$ & $+8.08716\cdot 10^{-10}$ \\ \hline
$3\cdot 10^{6}$ & $-8.22419\cdot 10^{-10}$ \\ \hline
$4\cdot 10^{6}$ & $+8.01923\cdot 10^{-10}$ \\ \hline
$5\cdot 10^{6}$ & $+2.78245\cdot 10^{-10}$ \\ \hline
$6\cdot 10^{6}$ & $-5.00102\cdot 10^{-10}$ \\ \hline
$7\cdot 10^{6}$ & $-5.21564\cdot 10^{-10}$ \\ \hline
\end{tabular}%
\bigskip
\end{center}

The residues in (\ref{osc}) may be expressed by complex zeros of $\zeta $ in
a manifest way. Introduce function $\xi (s)=\xi (1-s)$ as usual%
\begin{equation*}
\xi (s):=2(s-1)\pi ^{-s/2}\Gamma \left( 1+\frac{s}{2}\right) \zeta (s).
\end{equation*}%
Since $\xi $ may be factorized using Hadamard product%
\begin{equation*}
\xi (s)=\dprod\limits_{\rho }\left( 1-\frac{s}{\rho }\right)
\end{equation*}%
(the product being taken over paired complex zeros $\rho $) one can see that
the residue of $1/\xi (s)$ in each particular pole $\rho $ may be expressed
as a function of all remaining $\rho $s. Let us label consecutive complex
zeros as $\rho _{i}$ according to the increase of $\left\vert \gamma
\right\vert $. Choose further a particular zero $\rho _{n}$ and introduce
certain "crippled" function $\xi _{k}$ as%
\begin{eqnarray}
\xi (s) &=&\left( s-\rho _{n}\right) \xi _{n}(s)  \label{ksi} \\
\xi _{n}(s) &:&=-\frac{1}{\rho _{n}}\left( 1-\frac{s}{\rho _{n}^{\ast }}%
\right) \dprod\limits_{i=1}^{n-1}\left( 1-\frac{s}{\rho }\right)
\dprod\limits_{i=n+1}^{\infty }\left( 1-\frac{s}{\rho }\right)  \notag
\end{eqnarray}%
$1/\xi _{k}(s)$ is obviously regular at $\rho _{k}$ and%
\begin{equation}
\func{Res}\left( \frac{1}{\zeta (2s+2)};s=-\frac{1+\rho _{n}}{2}\right) =-%
\frac{1}{2}\frac{h\left( \rho _{n}^{\ast }\right) }{\xi _{n}(\rho _{n})}
\label{res}
\end{equation}%
where asterisk denotes complex conjugation and%
\begin{equation*}
h(s):=2(s-1)\pi ^{-s/2}\Gamma \left( 1+\frac{s}{2}\right)
\end{equation*}%
which is regular at $\rho $. Formula (\ref{res}) is not very useful in
practice since the product over paired complex zeros converges slowly,
nevertheless, when combined with (\ref{osc}), gives an explicit formula
relating any $\overset{\sim }{c}_{k}$ to complex zeros $\rho $.

Using different approach B\'{a}ez-Duarte gave another version of explicit
formula of this kind, cf. \cite{Duarte1}, Theorem 1.5. The present approach
is based entirely on N\"{o}rlund-Rice's integral and provides, I believe,
immediate and more natural interpretation for the existence of oscillatory
component of the coefficients $c_{k}$ which is completely hidden in their
primary definition (\ref{ck}). After all, the behavior of this very
component is crucial for the RH, see Fig. 5. Last not least, resulting
formulas are more handy and more effective in numerical calculations.

\section{Discussion and open questions}

Having any criterion for the RH the key thing is to say how "useful" it may
be in numerical experiments. For example, Li's criterion states that RH is
true if certain numbers $\lambda _{n}$ are positive \cite{Maslanka2}.
However, the fact that $n$ initial $\lambda _{n}$ are positive implies that
roughly only $\sqrt{n}$ complex zeros lay on the critical line \cite%
{Oesterle}. Since high $\lambda _{n}$ are extremely difficult to compute we
can honestly say that Li's criterion, although very elegant, is pretty
useless in practice. The natural question arises now: does the fact that $k$
initial $c_{k}$ obeys (\ref{criterion}) implies that certain number $m$ of
initial complex zeros lay on the critical line? What is the relation between 
$k$ and $m$? Below we present simple argument based on formulas derived in
this paper.

The consecutive residues (\ref{res}) seem to be bounded sequence of complex
numbers, cf. Fig. 6, so we neglect their influence. Using Euler beta
function expansion for large values of its first argument \cite{Wolfram} 
\begin{equation}
B(a,b)\propto \Gamma (b)a^{-b}\left[ 1-\frac{b(b-1)}{2a}\left( 1+O\left( 
\frac{1}{a}\right) \right) \right] ,\quad \left\vert a\right\vert
\rightarrow \infty  \label{asympt}
\end{equation}

we get:%
\begin{gather}
B\left( k+1,\frac{1+\rho }{2}\right) \propto (k+1)^{-\frac{1+\rho }{2}%
}\Gamma \left( \frac{1+\rho }{2}\right) =  \label{asympt1} \\
=\Gamma \left( \frac{1+\rho }{2}\right) (k+1)^{-\frac{3}{4}}\left[ \cos
\left( \frac{\gamma }{2}\log \left( k+1\right) \right) +i\sin \left( \frac{%
\gamma }{2}\log \left( k+1\right) \right) \right]  \label{asympt2}
\end{gather}%
It is clear that if $\gamma $ is real then we get bounded oscillations in $%
\overset{\sim }{c}_{k}$ and, in the view of \ (\ref{criterion}), RH is
satisfied. Suppose, however, that there exists somewhere defiant extremely
high zero $\frac{1}{2}+i\gamma $ with $\gamma $ having non-vanishing
imaginary part. (More precisely, there should be four such zeros at once
placed symmetrically with respect to the critical axis and the real axis.)
We are now sure that such a zero, if any, lies higher than $\gamma \simeq
10^{13}$ \cite{Gourdon}, \cite{Wedeniwski} and probably higher than $\gamma
\simeq 10^{21}$ \cite{Odlyzko}. This non-vanishing imaginary part would
cause the trigonometric functions in (\ref{asympt2}) to acquire, roughly
speaking, growing amplitudes thus violating criterion (\ref{criterion}).
Now, the amplitude of such growing oscillations would be extremely small due
to elementary properties of the Euler beta function which appears in (\ref%
{osc}). Specifically, for $\gamma =$ $2\cdot 10^{21}$ (which is not much
higher than the range acquired recently by Odlyzko) the amplitude due to the
corresponding zero is roughly $10^{-3.5\cdot 10^{21}}$. Of course, the task
to extract this particular harmonic, let alone to check whether its
amplitude grows, is far beyond any numerical capabilities. In other words,
even if the relatively low complex zero of zeta, say the hundredth, had some
small shift off the critical line we would not be capable of finding this
using (\ref{criterion}) since its amplitude, in the sense of (\ref{asympt2}%
), is about $10^{78}$ times smaller than the amplitude of the first zero.

\begin{center}
*
\end{center}

This is rather sad news for all those who wish to (dis)prove the celebrated
RH using simple criteria. On the other hand, estimations like that presented
above enable one to understand why this problem is so tough. The Holy Grail
and immortality it would bring to persistent searcher remain still out of
reach... Well, perhaps there exists even simpler explanation to this
mystery, explanation not involving all that terrible mathematics. \textit{%
All things work together for good to them that love God} says the
inscription on Bernhard Riemann's tomb taken from St. Paul's \textit{Epistle
to the Romans} (8:28). Shouldn't these words be taken for a clue? Yet, for a
century and a half many mathematicians would sell their souls for the proof.
Clearly, there is an apparent contradiction between such a deal and St.
Paul's words.

\bigskip 

\textbf{Figure captions}

Fig. 1. Coefficients $c_{k}$ for $k$ up to $10^{4}$ exhibit no apparent sign
of oscillations. After calculating numerically the first $10^{3}$
coefficients $c_{k}$ this looked as "a pleasant smooth curve" \cite{Duarte2}.

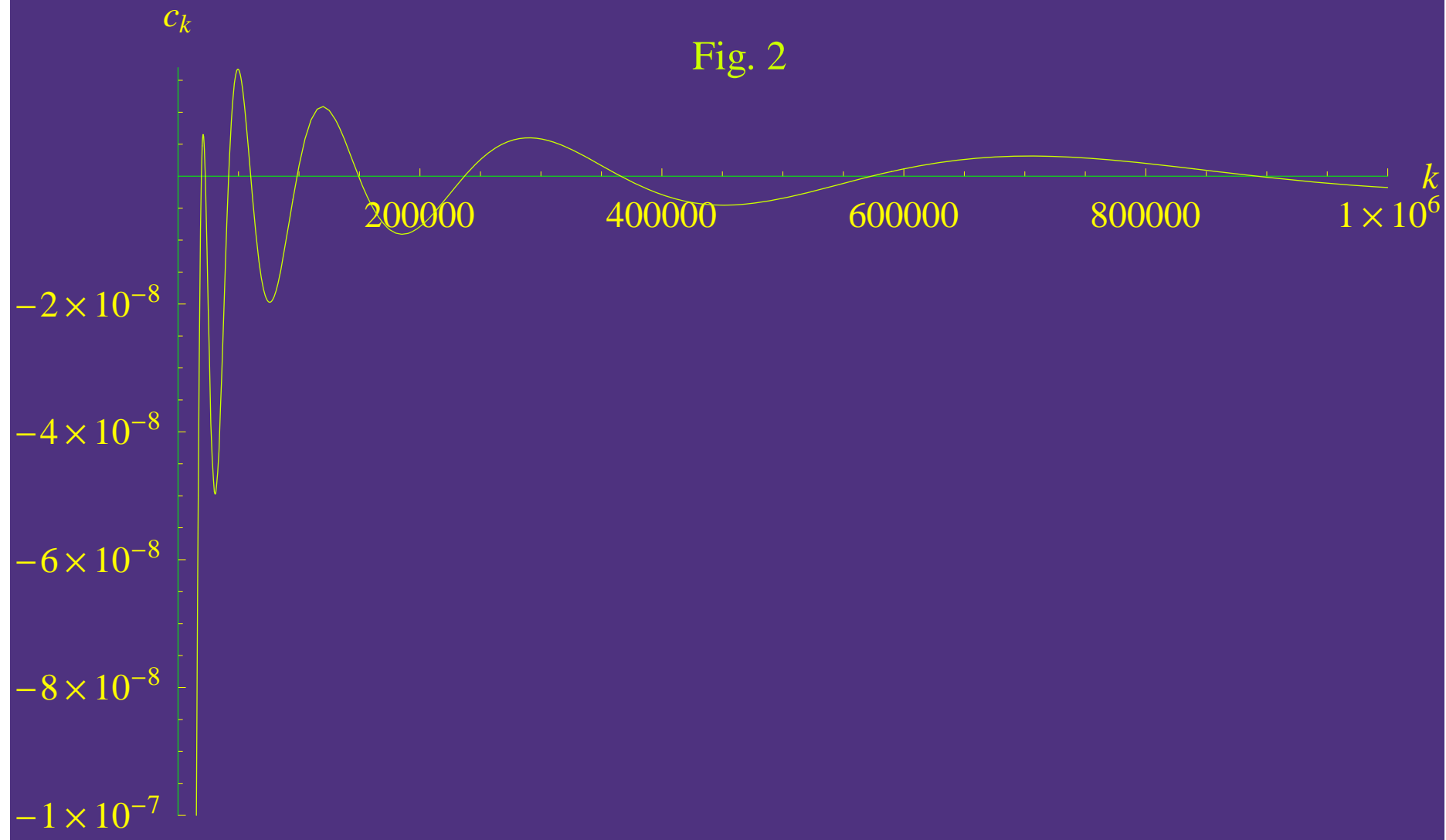
Fig. 2. Million coefficients $c_{k}$ computed using (\ref{ctrend}) and (\ref%
{osc}). Only 6 initial terms of the series were needed to accomplish
sufficient accuracy.

Fig. 3. The same as Fig. 2 with logarithmic scale in $k$.

Fig. 4. Behavior of $c_{k}$ (yellow) may be split into strictly growing
trend (red) plus the oscillating part (blue) here plotted in the logarithmic 
$k$-scale. Initially the trend dominates over the oscillations but it tends
to zero faster than these oscillations.

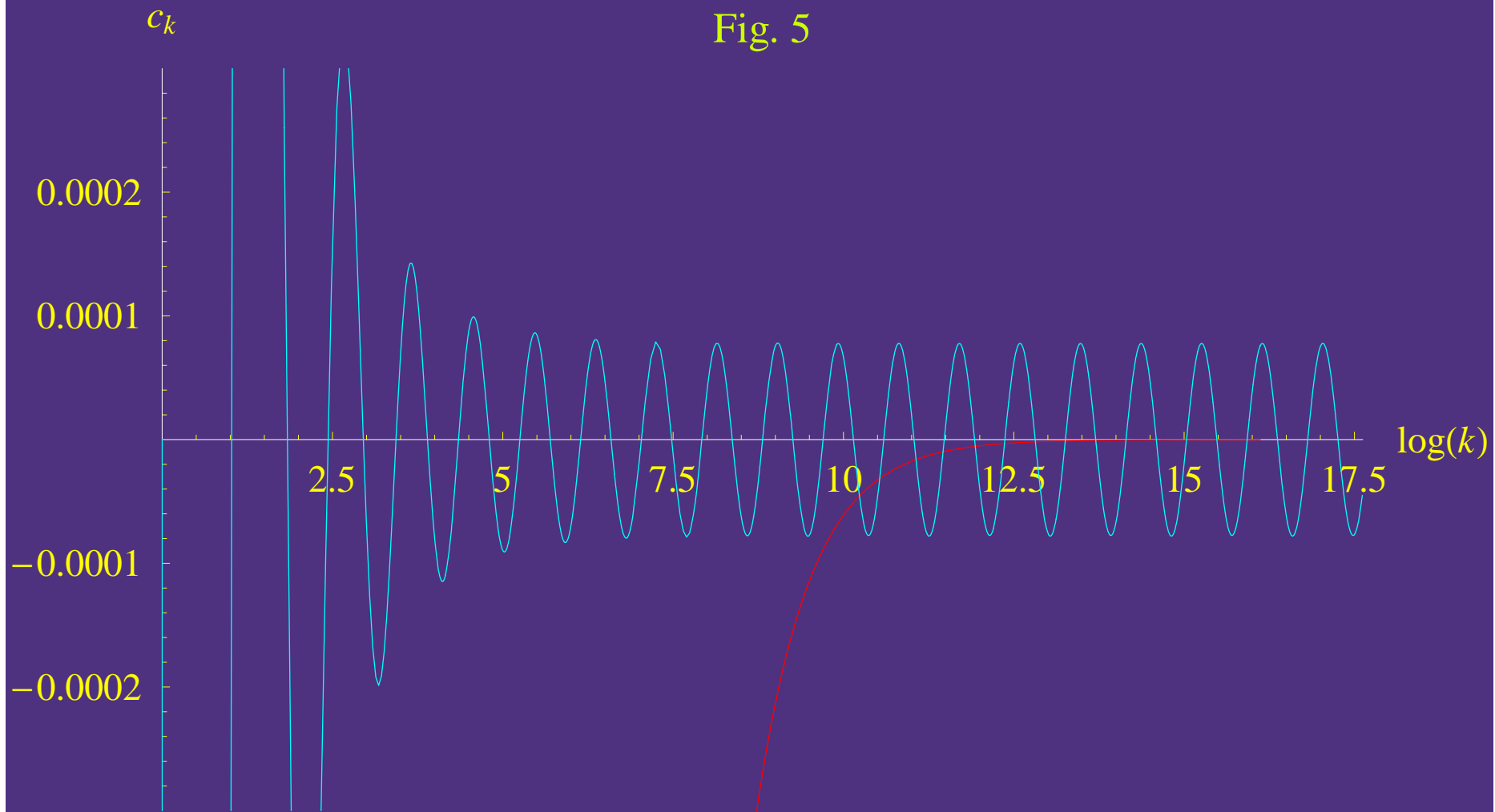
Fig. 5. Over 44 million components of $c_{k}$: trend $\overset{\_}{c}_{k}$
(red) and oscillating part $\overset{\sim }{c}_{k}$ (blue), both multiplied
by $k^{3/4}$ (compare RH criterion (\ref{criterion})). The former grows
strictly to zero whereas the latter seems to tend to perfect sine wave which
is dominated by the first zero. (Amplitudes of further zeros diminish
quickly.) If this tendency persisted to infinity \textit{for all zeros} RH
would be true. Maximum value of $k=44700000$ is related to parameter \texttt{%
\$MinNumber} in a particular version of \textit{Mathematica}.

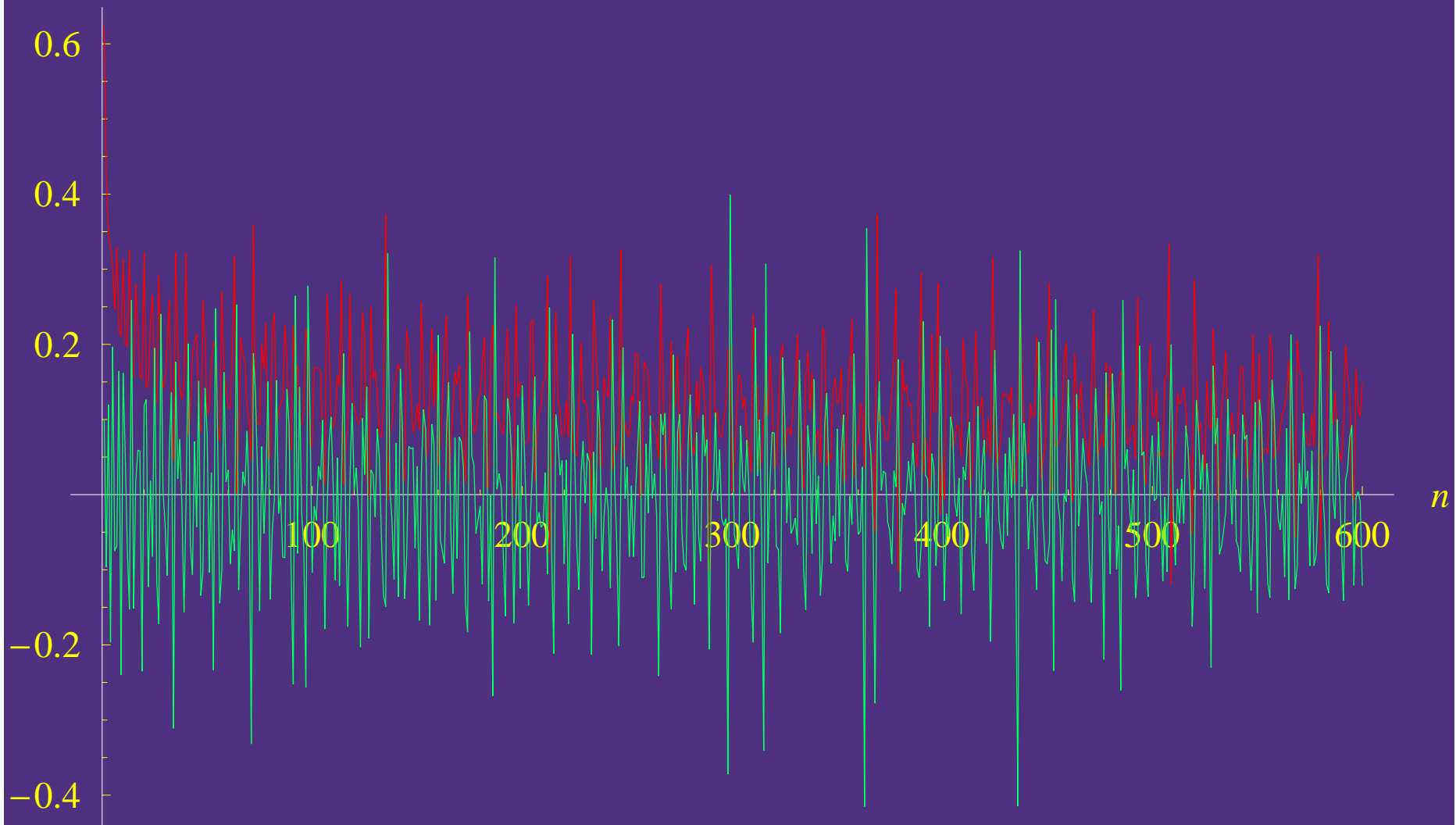
Fig. 6. 600 initial residues (\ref{res}). Red color denotes real part, green
-- imaginary part. Points are joined together in order to better visualize
their behavior.

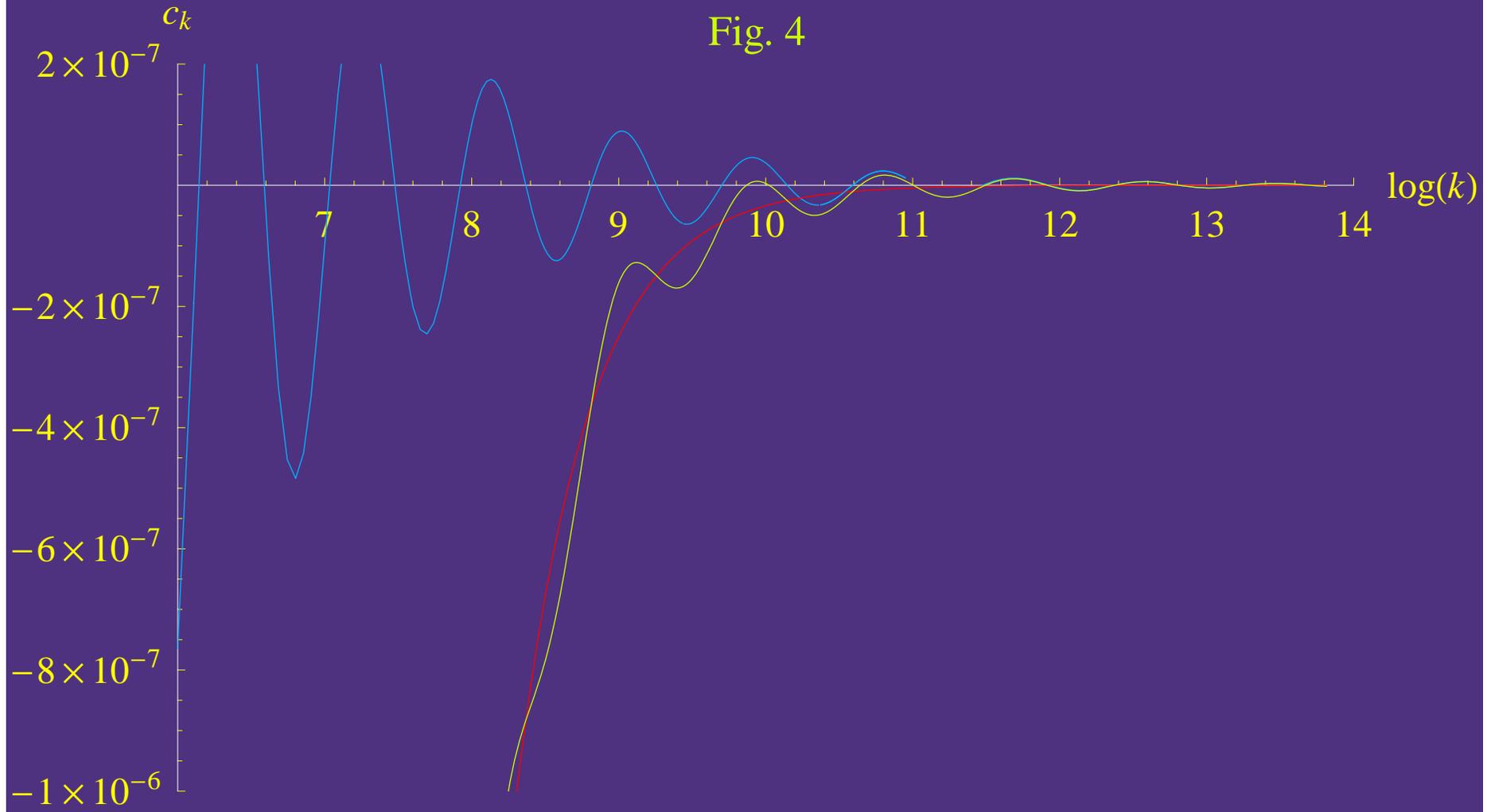
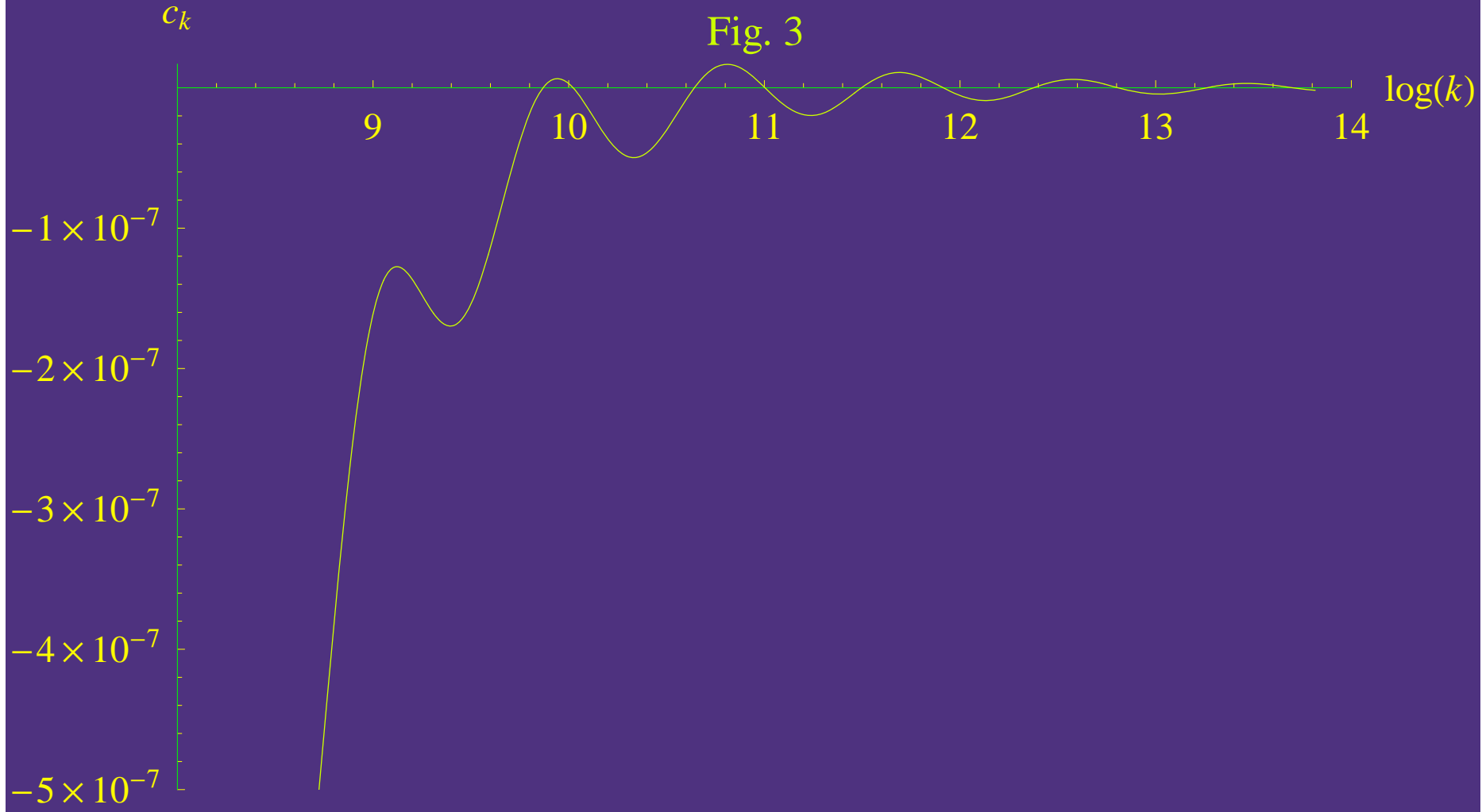
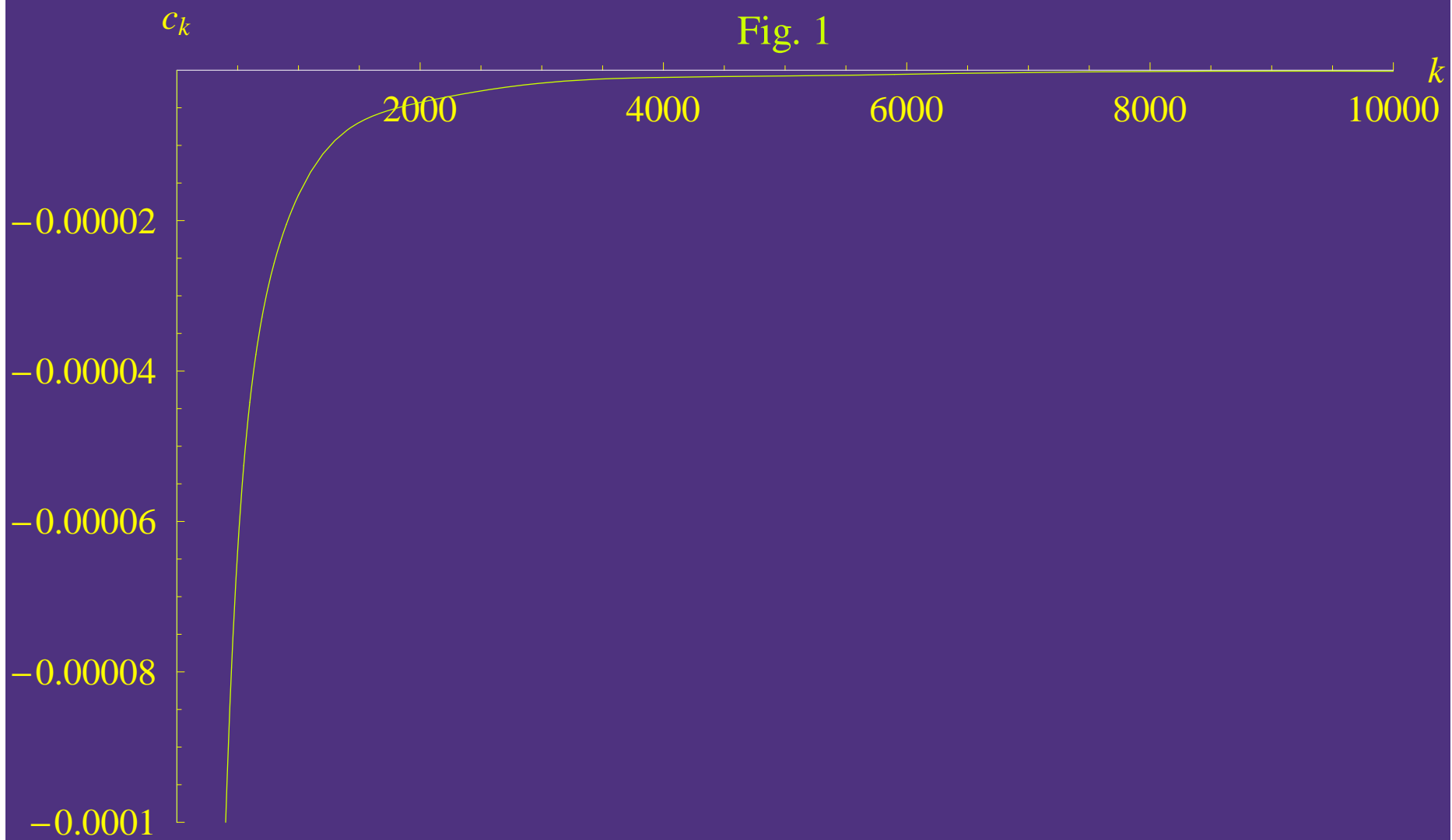
\end{document}